\documentclass[11pt]{amsart}
\usepackage{amssymb, amsmath, verbatim}
\theoremstyle{plain}
\newtheorem{thm}{Theorem}[section] 
\newtheorem{cor}[thm]{Corollary}
\newtheorem{prop}[thm]{Proposition}
\newtheorem{lem}[thm]{Lemma}
\theoremstyle{definition}

\theoremstyle{remark}
\newtheorem{rem}[thm]{Remark}

\numberwithin{equation}{section}%assigns eqns section numbers
\newcommand{\ind}{\operatorname{ind}}

\def\<{\left<}
\def\>{\right>}
\def\cstar{$C^*$-algebra}
\begin{document}
\title[$p$-summable commutators]{$p$-summable commutators in dimension $d$}
\author{William Arveson}
\thanks{supported by 
NSF grant DMS-0100487} 
\address{Department of Mathematics,
University of California, Berkeley, CA 94720}
\email{arveson@math.berkeley.edu}
%\date{26 March, 2003}
%
\begin{abstract} 
We show that many invariant subspaces $M$ for 
$d$-shifts $(S_1,\dots,S_d)$ of finite rank 
have the property that the orthogonal projection 
$P_M$ onto $M$ satisfies 
$$
P_MS_k-S_kP_M\in\mathcal  L^p,\qquad 1\leq k\leq d
$$
for every $p>2d$, $\mathcal  L^p$ denoting the 
Schatten-von Neumann class of all compact operators 
having $p$-summable singular value lists.  
In such cases, the $d$ tuple of operators 
$\bar T=(T_1,\dots,T_d)$ obtained by compressing 
$(S_1,\dots,S_d)$ to $M^\perp$ generates a 
$*$-algebra whose 
commutator ideal is contained in 
$\mathcal  L^p$ for every $p>d$. 

It follows that 
the $C^*$-algebra generated by $\{T_1,\dots,T_d\}$ 
and the identity is commutative modulo compact operators, 
the Dirac operator associated with 
$\bar T$ is Fredholm, 
and the index formula for the curvature invariant 
is stable under  compact perturbations and homotopy 
for this restricted class of finite rank $d$-contractions.   
Though this class is limited, 
we conjecture that the same conclusions persist under 
much more general circumstances.  
\end{abstract}
\maketitle

%\renewcommand{\baselinestretch}{1.3}\small
%use above for doublespacing 

\section{Introduction}\label{S:intro}

The purpose of this paper is to 
establish a result in higher dimensional operator theory 
that supports a general conjecture about the stability of 
the curvature invariant
under compact perturbations and homotopy.  
Specifically, we show that certain 
finite rank pure $d$-contractions 
$\bar T=(T_1,\dots,T_d)$ have the property that  
the \cstar\ generated by $\{T_1,\dots,T_n\}$ 
is commutative modulo compact operators.  
It follows that the Dirac operator 
associated with such a $d$-contraction is Fredholm, a key fact 
that leads to the desired stability properties for the 
curvature invariant by way of an index formula,  
see (\ref{Eq1}) below.  
These results represent a first step 
toward developing an effective Fredholm theory 
of $d$-contractions, an ingredient necessary for 
completing the index theorem that was partially 
established in \cite{arvDirac}, 
following up on \cite{arvCurv}.  

We first describe the issues 
that prompted this work.  
We use the term {\em multioperator} (of 
complex dimension $d=1,2,\dots$) to 
denote a $d$-tuple $\bar T=(T_1,\dots,T_d)$ of 
mutually commuting bounded operators acting on a 
common Hilbert space $H$.  Every multioperator 
$\bar T$ gives rise to an associated Dirac operator, 
whose definition we recall for the reader's convenience
(see  \cite{arvDirac} for more detail).  
Let $Z$ be a complex Hilbert space of 
dimension $d$ and let $\Lambda Z$ be the exterior 
algebra of $Z$, namely the direct sum of finite-dimensional 
Hilbert spaces 
$$
\Lambda Z = \sum_{k=0}^d \Lambda^k Z, 
$$ 
where $\Lambda^k Z$ denotes the $k$th exterior power of 
$Z$, and where $\Lambda^0 Z$ is defined to be 
the one-dimensional Hilbert space $\mathbb C$.  
For $k=1,\dots,d$, 
$\Lambda^k Z$ is spanned by wedge products of the 
form $z_1\wedge z_2\wedge\cdots \wedge z_k$, 
$z_i\in Z$, and 
the inner product in $\Lambda^k Z$ is 
uniquely determined by the formula 
$$
\langle z_1\wedge z_2\wedge\cdots \wedge z_k, 
w_1\wedge w_2\wedge\cdots \wedge w_k\rangle = \det (\langle z_i,w_j\rangle), 
$$
the right side denoting the determinant of the 
$k\times k$ matrix of inner products $\langle z_i,w_j\rangle$, 
$1\leq i,j\leq k$.  If we choose an orthonormal basis 
$e_1,\dots,e_d$ for $Z$, then there is a sequence of 
creation operators 
$C_1,\dots, C_d$ on $\Lambda Z$ 
that are defined uniquely by their action on 
generating vectors as follows
$$
C_i: z_1\wedge\cdots \wedge z_k\mapsto e_i\wedge z_1\cdots\wedge z_k,
$$
$z_1,\dots,z_k\in Z$, $1\leq k\leq d$, where each 
$C_i$ maps $\lambda\in\Lambda^0Z=\mathbb C$ 
to $\lambda e_i$ and 
maps the last summand $\Lambda^dZ$ to $\{0\}$.  
The operators $C_1,\dots,C_d$ satisfy the canonical 
anticommutation relations 
$$
C_iC_j+C_jC_i=0, \qquad C_i^*C_j+C_jC_i^*=\delta_{ij}\mathbf 1.  
$$

The Dirac operator of $\bar T$ 
is a self-adjoint operator $D$ acting on the Hilbert space  
$\tilde H=H\otimes \Lambda Z$ as follows: $D=B+B^*$, where 
$B$ is the sum 
$$
B=T_1\otimes C_1+T_2\otimes C_2+\cdots +T_d\otimes C_d.  
$$
If one replaces $e_1,\dots,e_d$ with 
a different orthonormal basis for $Z$, then of course one changes 
$D$; but the two Dirac operators are 
naturally isomorphic in a sense that we will not 
spell out here (see \cite{arvDirac}).  Thus the 
Dirac operator of $\bar T$ is uniquely determined 
by $\bar T$ up to isomorphism.  

The multioperator $\bar T$ is said to be {\em Fredholm} 
if its Dirac operator $D$ is a Fredholm operator.  
Since $D$ is self-adjoint, this simply means that it 
has closed range and finite-dimensional kernel.  In this 
case there is an integer invariant associated with 
$\bar T$, called the {\em index}, that is defined as 
follows.  Consider the natural $\mathbb Z_2$-grading 
of $\tilde H$, defined by the orthogonal decomposition 
$\tilde H=\tilde H_+\oplus \tilde H_-$, where 
$$
\tilde H_+=\sum _{k{\rm\ even}}H\otimes \Lambda^k Z,\qquad 
\tilde H_-=\sum _{k{\rm\ odd}}H\otimes \Lambda^k Z.  
$$
One finds that $D$ is an odd operator relative to this 
grading in the sense that $D\tilde H_+\subseteq \tilde H_-$, 
$D\tilde H_-\subseteq \tilde H_+$.   
Thus the decomposition 
$\tilde H=\tilde H_+\oplus \tilde H_-$ gives 
rise to a $2\times 2$ matrix representation 
$$
D=
\begin{pmatrix}
0&D_+^*\\
D_+&0
\end{pmatrix}, 
$$
$D_+$ denoting the restriction of $D$ to $\tilde H_+$.  
When  $D$ is 
Fredholm, one finds that $D_+\tilde H_+$ is a closed subspace 
of $\tilde H_-$ of finite codimension and 
$D_+$ has finite 
dimensional kernel.  Indeed, $\bar T$ is 
a Fredholm multioperator iff $D_+$ is a 
Fredholm operator in $\mathcal B(\tilde H_+,\tilde H_-)$.  
The index of $D_+$, namely 
$$
\ind (D_+)=\dim (\ker D\cap \tilde H_+) -\dim (\tilde H_-/D\tilde H_+), 
$$
is an integer invariant for 
Fredholm multioperators that is stable under compact 
perturbations and homotopy.  

By analogy with the index theorems of Atiyah and Singer,  
one might expect that 
the computation of the index of Fredholm multioperators 
will lead to important relations between the geometric 
and analytic properties of multioperators, and perhaps 
connect with basic issues of algebraic geometry.  
Such a program requires that one 
should have effective tools for a) determining when a given 
Dirac operator is Fredholm and b) computing the index 
in terms of concrete geometric properties of 
its underlying multioperator.  
Some progress has been made in the direction of b), and 
we will describe that in Remark \ref{rem1} below.  
However, the problem a) of 
proving that the natural examples of multioperators are 
Fredholm remains largely open.  It is that problem 
we want to address in this paper.

Here is a useful sufficient 
criterion for Fredholmness; we reiterate the proof given in 
\cite{arvLuminy} for the reader's convenience.  

\begin{prop}\label{prop1}
Let $\bar T=(T_1,\dots, T_d)$ be a multioperator 
satisfying    
\begin{enumerate}
\item[(i)]$\bar T$ is essentially normal in that 
all self-commutators 
$T_kT_j^*-T_j^*T_k$ are compact, $1\leq j, k\leq d$, and   
\item[(ii)]$T_1T_1^*+\cdots+T_dT_d^*$ is a Fredholm operator. 
\end{enumerate}
Then $\bar T$ is a Fredholm multioperator.  
\end{prop}

\begin{proof}
Since the Dirac operator $D$ of $\bar T$ 
is self-adjoint, it suffices to 
show that $D^2$ is a Fredholm operator.  
To that end, consider 
$B=T_1\otimes C_1+\cdots+T_d\otimes C_d$.  Since 
$T_j$ commutes with $T_k$ and $C_j$ anticommutes with 
$C_k$, a straightforward computation shows that 
$B^2=0$.  Hence 
$$
D^2=(B+B^*)^2=B^*B+BB^* = 
\sum_{k,j=1}^d T_k^*T_j\otimes C_k^*C_j+
\sum_{k,j=1}^d T_jT_k^*\otimes C_jC_k^*.  
$$
Using 
$C_jC_k^*=\delta_{jk}\mathbf 1-C_k^*C_j$, we can write 
the second term on the right as 
$$
F\otimes\mathbf 1-\sum_{k,j=1}^d T_jT_k^*\otimes C_k^*C_j, 
$$
where $F=T_1T_1^*+\cdots+T_dT_d^*$, so that 
$$
D^2=F\otimes \mathbf 1 +\sum_{k,j=1}^d (T_k^*T_j-T_jT_k^*)\otimes C_k^*C_j.  
$$
Since $F\otimes \mathbf 1$ is a Fredholm operator 
by (ii)  and 
each summand in the second term is compact by (i), it follows 
that $D^2$ is a Fredholm operator.  
 \end{proof}

\begin{rem}[Finite rank $d$-contractions]\label{rem1}
We are primarily concerned with finite rank 
$d$-contractions, that is, 
multioperators $\bar T=(T_1,\dots, T_d)$ 
that define row contractions in the sense 
that $T_1T_1^*+\dots+T_dT_d^*\leq\mathbf 1$,  
whose defect operators 
$\mathbf 1-T_1T_1^*-\dots-T_dT_d^*$ have finite rank.  
$\bar T$ is said to be {\em pure} if the powers 
of the completely positive map $\phi(X)=\sum_k T_kXT_k^*$ 
satisfy $\phi^n(\mathbf 1)\downarrow 0$ as $n\to \infty$.  

Proposition \ref{prop1} implies that the Dirac 
operator of a finite rank $d$-contraction is Fredholm 
provided that the self-commutators $T_kT_j^*-T_j^*T_k$ are 
compact for all $1\leq k,j\leq d$.  
In that case, the $C^*$-algebra $C^*(T_1,\dots,T_d)$ generated by 
$\{T_1,\dots,T_d\}$ and the identity operator is commutative modulo  
compact operators $\mathcal K\subseteq \mathcal B(H)$, and we have an exact 
sequence of $C^*$-algebras
\begin{equation}\label{Eq13}
0\longrightarrow 
\mathcal K\longrightarrow C^*(T_1,\dots,T_d)+\mathcal K \longrightarrow 
C(X)\longrightarrow 0, 
\end{equation}
$X$ being a compact subset of the unit $(2d-1)$-sphere in $\mathbb C^d$.  

For every finite rank $d$-contraction 
$\bar T$ it is possible to define a 
real number $K(\bar T)$ in the interval 
$[0,{\rm rank\;}\bar T]$, called the 
curvature invariant.  $K(\bar T)$ is 
a geometric invariant of $\bar T$, 
defined as the integral 
of the trace of a certain matrix-valued function over the unit 
sphere in $\mathbb C^d$ (see \cite{arvCurv} for more detail).  
$K(\bar T)$ was computed for many examples in \cite{arvCurv}, 
and it was found to be an integer, namely the Euler characteristic 
of a certain finitely generated module over the polynomial 
ring $\mathbb C[z_1,\dots,z_d]$.  However, it was also shown 
in \cite{arvCurv} that this formula equating 
$K(\bar T)$ to an Euler characteristic fails to hold in 
general.  Still, that formula 
provided enough evidence to lead us to conjecture 
that $K(\bar T)$ is an integer in general, and that 
has now been established by Greene, Richter and 
Sundberg in \cite{GreRiSun}.  Unfortunately, 
the integer arising in the 
proof of the latter result - 
namely the rank of an almost-everywhere 
constant rank projection defined on the unit sphere 
of $\mathbb C^d$ - appeared to have no direct 
connection with spatial properties of the underlying 
$d$-tuple of operators $T_1,\dots,T_d$.  What was still lacking 
was a formula that relates $K(\bar T)$ to some natural integer 
invariant of $\bar T$ that holds in general and 
is, hopefully, easy to compute.  

Such considerations led us to initiate a new approach 
in \cite{arvDirac}.  Our aim was to 
introduce an appropriate notion 
of Dirac operator 
and to seek a formula that would relate 
the curvature invariant to the index 
of the associated Dirac operator - hopefully in general. 
Assuming that 
$\bar T$ is a {\em graded} $d$-contraction, 
then the result of \cite{arvDirac} 
is that both $\ker D_+$ and $\ker D_+^*$ are 
finite-dimensional, and moreover 
\begin{equation}\label{Eq1}
(-1)^d K(\bar T)=\dim\ker D_+ - \dim\ker D_+^*.  
\end{equation}
Significantly, there are no known exceptions 
to this formula; in particular, it persists 
for the (ungraded) examples that violated 
the previous formula that related $K(\bar T)$ 
to an Euler characteristic.  
As we have described in \cite{arvDirac}, 
it is natural to view 
(\ref{Eq1}) as an operator-theoretic 
counterpart of the Gauss-Bonnet-Chern 
formula of Riemannian geometry in its 
modern dress as an index theorem 
(see page 311 of \cite{GilMur}).

Notice, however, that even in the 
graded cases the right 
side of (\ref{Eq1}) is unstable 
if $D$ is not a Fredholm operator, 
despite the fact that both 
subspaces $\ker D_+$ and $\ker D_+^*$ must be 
finite-dimensional.  
On the other hand, if $D$ is a Fredholm operator 
then (\ref{Eq1}) reduces to a stable formula 
\begin{equation}\label{Eq2}
(-1)^d K(\bar T)=\ind D_+.  
\end{equation}
In view of Proposition \ref{prop1} and 
the stability properties of the index of 
Fredholm operators one 
may conclude: {\em Within the class of finite rank 
graded $d$-contractions $\bar T$ whose self-commutators 
are compact, the curvature 
invariant $K(\bar T)$ is stable under compact perturbations 
and homotopy.}
\end{rem}

Unfortunately, it is not known if the self-commutators of 
pure finite-rank graded $d$-contractions are always compact, 
though we believe that they are.  
More generally, we believe that 
the Dirac operator of any 
pure finite rank $d$-contraction $\bar T$ -- graded 
or not -- is a Fredholm operator and, moreover, 
that formula (\ref{Eq2}) continues 
to hold in that generality.  
Precise formulations of these conjectures 
will be found in Section \ref{S:conc}. 
We will prove the most tractable cases of the first 
of these two conjectures in Sections \ref{S:statement} 
and \ref{S:proofs}.

\section{Statement of Results}\label{S:statement}

Let $\bar S=(S_1,\dots,S_d)$ be the $d$-shift of 
rank one, i.e., the multioperator that 
acts on the symmetric Fock 
space $H^2(\mathbb C^d)$ over the $d$-dimensional one-particle 
space $\mathbb C^d$ by symmetric tensoring with a fixed 
orthonormal basis $e_1,\dots,e_d$ for $\mathbb C^d$.  
We lighten notation by writing $H^2$ for $H^2(\mathbb C^d)$, 
the dimension $d$ being 
a positive integer (normally larger than $1$) that will be fixed 
throughout.  The elements of $H^2$ can be realized as certain 
holomorphic functions defined in the open unit ball 
$$
B_d=\{z=(z_1,\dots,z_d)\in\mathbb C^d: 
|z|=(|z_1|^2+\dots+|z_d|^2)^{1/2}<1\},
$$
and in this function-theoretic realization 
the rank one $d$-shift is the $d$-tuple of multiplication 
operators $S_k:f(z)\mapsto z_kf(z)$, $1\leq k\leq d$ 
(see \cite{arvSubalgIII}).

Let $r$ be a positive integer, let $E$ be a Hilbert space of 
dimension $r$, and consider the  
$d$-tuple of operators defined defined on 
$H^2\otimes E$ by 
$$
S_k\otimes\mathbf 1_E:
f\otimes\zeta\mapsto z_kf\otimes \zeta, \qquad 1\leq k\leq d.  
$$
It will be convenient to overwork notation by writing  
$\bar S=(S_1,\dots, S_d)$ for these operators as well, and to 
refer to that multioperator as the $d$-shift 
of rank $r$.  
Thus the $d$-shift $(S_1,\dots,S_d)$ of rank $r$ acts as follows 
on elements of the form $f\otimes\zeta$, with $f\in H^2$ 
and $\zeta\in E$:
$$
S_k: f\otimes\zeta\mapsto z_kf\otimes\zeta, \qquad 1\leq k\leq d.  
$$
The $d$-shift is known to be universal in the class 
of pure $d$-contractions in the sense that {\em every 
pure $d$-contraction $\bar T=(T_1,\dots,T_d)$ of rank $r$ 
is unitarily equivalent 
to one obtained by compressing the 
$d$-shift of rank $r$ to the orthogonal complement of an 
invariant subspace} (see \cite{arvSubalgIII}).

We will make use of the natural partial ordering on the discrete abelian 
group $\mathbb Z^d$; for $m=(m_1,\dots,m_d), n=(n_1,\dots,n_d)\in\mathbb Z^d$
we write $m\leq n$ if $m_k\leq n_k$ for every $k=1,\dots,d$.  For every 
$n\geq 0$ in $\mathbb Z^d$ there is a monomial in $H^2$ 
defined by
$$
z^n=z_1^{n_1}z_2^{n_2}\dots z_d^{n_d}. 
$$
The set of all monomials $\{z^n: n\geq 0\}$ form an orthogonal (but not orthonormal) 
set which spans $H^2$.  Similarly, an element $\xi\in H^2\otimes E$ having 
the particular form
$
\xi = z^n\otimes\zeta 
$
where $n\geq 0$ and $\zeta\in E$ is called a {\em monomial} in $H^2\otimes E$.  
Notice that monomials of the form $z^m\otimes\eta$ and $z^n\otimes\zeta$ are 
orthogonal if $m\neq n$, or if $m=n$ and $\eta\perp\zeta$.  Obviously, 
$H^2\otimes E$ is spanned by monomials.  

\begin{thm}\label{oldThmA}
Let $M\subseteq H^2\otimes \mathbb C^r$ be an invariant subspace for 
the $d$-shift of rank $r$ that is generated as an invariant subspace by any 
set of monomials, and let $P_M$ be the orthogonal projection onto $M$.  
Then for every $p>2d$, 
$$
P_MS_k-S_kP_M\in \mathcal L^p, \qquad 1\leq k\leq d.  
$$
\end{thm}

\begin{cor}\label{oldCor}
Let $M\subseteq H^2\otimes \mathbb C^r$ be an invariant subspace
satisfying the hypotheses of Theorem \ref{oldThmA}, let 
$\bar T=(T_1,\dots,T_d)$ be the $d$-contraction obtained by compressing 
$\bar S$ to the Hilbert space $(H^2\otimes E)/M\cong M^\perp$, and let 
$\mathcal A$ be the $*$-algebra generated by $\{T_1,\dots,T_d\}$ 
and the identity.  Then 
\begin{enumerate}
\item[(i)] the commutator ideal of $\mathcal A$ is 
contained in $\mathcal L^p$ for every 
$p>d$, and  
\item[(ii)] $\bar T$ is a Fredholm multioperator.  
\end{enumerate}
\end{cor}

\begin{rem}[Structure of quotient modules]
Consider the two dimensional rank one case $d=2$, 
$E=\mathbb C$.  The simplest nontrivial 
example of an invariant subspace 
satisfying the hypotheses of Theorem \ref{oldThmA} is the subspace 
$M\subseteq H^2(\mathbb C^2)$ generated by the single monomial 
$f(z_1,z_2)=z_1z_2$.  In this case it is possible to compute the 
operators $T_1,T_2\in\mathcal B(H^2/M)$ 
in explicit terms.  Once that is done, 
one can verify directly that for this example the 
self-commutators $T_i^*T_j-T_jT_i^*$, $1\leq i,j\leq 2$, 
are in $\mathcal L^p$ for every $p>2$, so that Proposition 
\ref{prop1} implies that $(T_1,T_2)$ is a Fredholm pair.  

However, for more general examples the quotient modules
$H^2/M$ are not recognizable (consider 
the invariant subspace $M\subseteq H^2(\mathbb C^3)$ generated 
by the two monomials $z_1z_2$ and $z_2z_3^2$) .
Thus our proof of Theorem 
\ref{oldThmA}  and Corollary \ref{oldCor} is based on 
different ideas.  
\end{rem}

\begin{rem}[Structure of the algebraic set $X$]
Let $f_1,\dots, f_s$ be a set of homogeneous 
polynomials in $\mathbb C[z_1,\dots,z_d]$ and 
consider the closed 
invariant subspace $M\subseteq H^2(\mathbb C^d)$ 
that they generate.  Let $\bar T=(T_1,\dots,T_d)$ be the 
rank-one $d$-contraction associated with the quotient 
$H^2/M$.  Assuming that the 
self-commutators $T_j^*T_k-T_kT_j^*$ are compact, 
$1\leq j,k\leq d$, then we have an exact sequence 
of \cstar s (\ref{Eq13}) 
that terminates in $C(X)$, where $X$ 
is the projective algebraic 
set defined by the common zeros of $f_1,\dots,f_s$.   
However, if the $f_k$ are all 
monomials then $X$ is trivial - a union 
of coordinate axes.  Thus 
the Fredholm $d$-tuples provided by 
Corollary \ref{oldCor} fail to make significant connections 
with algebraic geometry.  
\end{rem}

\section{Proofs.}\label{S:proofs}

Turning now to the proof of Theorem \ref{oldThmA} 
and Corollary \ref{oldCor}, let $E=\mathbb C^r$ and let   
$M\subseteq H^2\otimes E$ be an invariant subspace 
generated by some set of monomials.  We have to show that 
the commutators 
$$
[P_M,S_k] = P_MS_k-S_kP_M, \qquad 1\leq k\leq d
$$
belong to $\mathcal L^p$ for $p>2d$.  
Because of the obvious symmetry we treat only the case 
$k=1$.  We will make repeated use of the following 
elementary property of the Schatten-von Neumann 
classes $\mathcal L^p$: For any 
pair of Hilbert spaces $H_1$, $H_2$, any operator 
$B\in\mathcal B(H_1,H_2)$, and any $p\geq 1$ 
one has 
\begin{equation}\label{Eq12}
B\in\mathcal L^{2p}\iff B^*B\in \mathcal L^p\iff 
BB^*\in\mathcal L^p.  
\end{equation}
Thus, in order to prove Theorem \ref{oldThmA} it suffices to show 
that the operator 
$$
B=[P_M,S_1]^*=S_1^*P_M-P_MS_1^*=(\mathbf 1-P_M)S_1^*P_M
$$
satisfies $B^*B\in \mathcal L^p$ for every $p>d$.  Equivalently, 
we will prove: 
\begin{thm}\label{oldAss2.1}
Let $A$ be the restriction of the operator $(\mathbf 1-P)S_1^*$ 
to $M$.  Then $A^*A\in\mathcal L^p$ for every $p>d$.  
\end{thm}
Consider the natural decomposition of $M$ induced by $S_1$
$$
M = \overline{S_1M}\oplus (M\ominus S_1M).  
$$
Since $S_1M\perp M\ominus S_1M$, we have 
\begin{equation}\label{oldEq2.7}
S_1^*(M\ominus S_1M)\subseteq M^\perp,   
\end{equation}
and therefore both $S_1^*$ and $(\mathbf 1-P_M)S_1^*$ restrict 
to the same operator on $M\ominus S_1M$.  
Thus it suffices to establish the following two results, the 
principal one being Lemma \ref{oldProp2.3}.  
\begin{lem}\label{oldProp2.2}
$S_1^*S_1$ leaves $M$ invariant, hence the restriction of 
$(\mathbf 1-P_M)S_1^*$ to $\overline{S_1M}$ is zero.  
\end{lem}

\begin{lem}\label{oldProp2.3}
The restriction of $S_1^*$ to $M\ominus S_1M$ belongs 
to $\mathcal L^{2p}$ for every $p>d$.  
\end{lem}

We first bring in 
an action of the compact group $\mathbb T^d$ that will be useful.  
The full unitary group of the one-particle space 
$\mathbb C^d$ acts naturally as unitary operators 
on $H^2=H^2(\mathbb C^d)$, and by restricting that  
representation to the abelian subgroup of all unitary operators 
which are diagonal relative to the usual orthonormal 
basis for $\mathbb C^d$, one obtains a strongly continuous
unitary representation of the $d$-dimensional torus 
$\Gamma_0:\mathbb T^d\mapsto \mathcal B(H^2)$.  In more explicit 
terms, if we consider the elements of 
$H^2(\mathbb C^d)$ as holomorphic functions 
defined on the open unit ball $B_d\subseteq\mathbb C^d$, 
the action of $\Gamma_0$ is given by 
$$
\Gamma_0(\lambda): f(z_1,\dots,z_d)\mapsto f(\lambda_1z_1,\dots,\lambda_d z_d),
\qquad \lambda=(\lambda_1,\dots,\lambda_d)\in\mathbb T^d.
%\tag{1.2}
$$
By increasing the multiplicity appropriately, we obtain a 
corresponding representation $\Gamma: \mathbb T^d\to\mathcal B(H^2\otimes E)$, 
\begin{equation}\label{old2.4}
\Gamma(\lambda)=\Gamma_0(\lambda)\otimes\mathbf 1_E,\qquad \lambda\in\mathbb T^d.  
\end{equation}
We have the following relations between $\Gamma$ and 
the rank $r$ $d$-shift
\begin{equation}\label{oldEq2.4}
\Gamma(\lambda)S_k\Gamma(\lambda)^*=\lambda_kS_k,\qquad k=1,\dots,d,\quad
\lambda=(\lambda_1,\dots,\lambda_d)\in\mathbb T^d.  
%\tag{2.4}
\end{equation}

The character group of $\mathbb T^d$ is the discrete 
abelian group $\mathbb Z^d$, an element $n=(n_1,\dots,n_d)$ in 
$\mathbb Z^d$ being associated with the following character of $\mathbb T^d$,
$$
\lambda\mapsto \lambda^n=\lambda_1^{n_1}\dots\lambda_d^{n_d}
\in\mathbb T, \qquad 
\lambda=(\lambda_1,\dots,\lambda_d)\in\mathbb T^d.  
$$
Notice that every monomial $z^n\otimes \zeta$, $n\geq 0$ is an eigenvector 
of $\Gamma$, 
$$
\Gamma(\lambda)(z^n\otimes\zeta) = \lambda^n\cdot z^n\otimes\zeta, 
\qquad \lambda\in\mathbb T^d.
$$
Conversely, if $\xi\in H^2\otimes E$ satisfies 
$\Gamma(\lambda)\xi=\lambda^n\xi$ for every $\lambda\in\mathbb T^d$ and 
if $\xi\neq 0$ then we must have $n\geq 0$ and $\xi$ must be a monomial of 
degree $n$.  

The role of $\Gamma$ is described in the following 
proposition.  

\begin{prop}\label{oldProp2.5}
For any closed subspace $M\subseteq H^2\otimes E$ that is invariant 
under $S_1,\dots,S_d$, the following are equivalent.
\begin{enumerate}
\item[(i)] 
$M$ is generated (as a closed $\{S_1,\dots,S_d\}$-invariant 
subspace) by a set of monomials in $H^2\otimes E$.  
\item[(ii)]
$\Gamma(\lambda)M\subseteq M$, for every $\lambda\in\mathbb T^d$.  
\end{enumerate}
\end{prop}

\begin{proof} (i)$\implies$(ii):  
Obviously, operators of the form 
$S_1^{n_1}S_2^{n_2}\dots S_d^{n_d}$, with 
$n_1,\dots,n_k\geq 0$
must map monomials in $H^2\otimes E$ to other monomials in $H^2\otimes E$, 
and since monomials in $H^2\otimes E$ form one-dimensional $\Gamma$-invariant 
subspaces, 
it follows that any invariant subspace $M$ that is generated by 
a set of monomials must also be invariant under the action of $\Gamma$.  

(ii)$\implies$(i): 
Any closed linear subspace $M\subseteq H^2\otimes E$ that is invariant 
under the action of $\Gamma$  must be spanned 
by its spectral subspaces
$$
M(n)=\{\xi\in M: \Gamma(\lambda)\xi=\lambda^n\xi,\quad\lambda\in\mathbb T^d\}, 
$$
for $n\in\mathbb Z^d$.  We have already pointed out that 
if such a subspace $M(n)$ is not $\{0\}$ then one must have $n\geq 0$, and 
that it must have the form 
$$
M(n)= z^n\otimes E_0=\{z^n\otimes\zeta: \zeta\in E_0\}, 
$$ 
$E_0$ being some subspace of $E$.  
Thus $M$ is spanned by the monomials it 
contains, and in particular it is generated as in (i).  
\end{proof}

\begin{proof}[proof of Lemma \ref{oldProp2.2}]
Proposition \ref{oldProp2.5} implies that $M$ is invariant under the von Neumann 
algebra generated by the range $\Gamma(\mathbb T^d)$ of $\Gamma$, and thus 
it suffices to show that $S_1^*S_1$ belongs to that algebra.  Because of the 
double commutant theorem it is enough to show that for every 
operator $T$ satisfying 
$$
\Gamma(\lambda)T=T\Gamma(\lambda), \qquad \lambda\in \mathbb T^d, 
$$ 
we have $TS_1^*S_1=S_1^*S_1T$.  

Now $H^2\otimes E$ decomposes into an orthogonal direct sum of spectral 
subspaces for $\Gamma$, namely the subspaces of the form $z^n\otimes E$ where 
$n\in\mathbb Z^d$ satisfies $n\geq 0$, and by virtue of its commutation relation 
with $\Gamma$, $T$ must leave each of these subspaces invariant.  
Thus there is a sequence of operators 
$T_n\in\mathcal B(E)$, $n\in\mathbb Z^d$, $n\geq 0$ such that the restriction of $T$ to 
$z^n\otimes E$ is given by 
\begin{equation}\label{oldEq2.6}
T(z^n\otimes \zeta)=z^n\otimes T_n\zeta, \qquad \zeta\in E.  %\tag{2.6}
\end{equation}

Consider now the action of $S_1^*S_1$ on 
the spectral subspace $z^n\otimes E$.  Writing 
$n=(n_1,\dots,n_d)$ with $n_k\geq 0$ we have 
$$
S_1^*S_1(z^n\otimes\zeta)=S_1^*(z_1^{n_1+1}z_2^{n_2}\dots z_d^{n_d}\otimes\zeta)=
(S_1^*(z_1^{n_1+1}z_2^{n_2}\dots z_d^{n_d}))\otimes\zeta.  
$$
Using formula (3.9) of \cite{arvSubalgIII} 
we have 
$$
S_1^*(z_1^{n_1+1}z_2^{n_2}\dots z_d^{n_d})=\frac{n_1+1}{|n|+1}\cdot 
z_1^{n_1}z_2^{n_2}\dots z_d^{n_d}, 
$$
where $|n|$ denotes $n_1+n_2+\dots+n_d$.  Thus 
$$
S_1^*S_1(z^n\otimes\zeta)=\frac{n_1+1}{|n|+1}\cdot (z^n\otimes\zeta).  
$$
Thus the restriction of $S_1^*S_1$ to each spectral subspace of 
$\Gamma$ is a scalar multiple of the identity; because of 
(\ref{oldEq2.6}), $S_1^*S_1$ must commute with $T$.
\end{proof}

\begin{proof}[proof of Lemma \ref{oldProp2.3}]
We prove Lemma \ref{oldProp2.3} in two 
assertions as follows.  
\vskip0.1in
{\bf Assertion 1.}
{\em 
There is a positive integer $q$ such that 
$$
M\ominus S_1M\subseteq
\overline{span}\{z_1^{n_1}z_2^{n_2}\dots z_d^{n_d}\otimes\zeta: 
(n_1,n_2,\dots,n_d)\geq 0,\  n_1\leq q,\ \zeta\in E\}.  
$$
}%
\vskip0.1in
{\bf Assertion 2.}
{\em For every positive integer $q$
the restriction $B$ of $S_1^*$ to 
$$
\overline{span}\{z^n\otimes\zeta: n\geq 0,\quad n_1\leq q,
\quad\zeta\in E\}
$$
satisfies $B^*B\in\mathcal L^p$ for every $p>d$.}  
\vskip0.1in
\begin{proof}[proof of Assertion 1]
We remark first that $M$ is finitely generated in the sense that 
there is a finite set $F$ of monomials in $M$ such that 
$$
M=\overline{\text span}\{f(S_1,\dots,S_d)\xi: \xi\in F, 
\quad f\in\mathbb C[z_1,\dots,z_d]\}.  
$$
Indeed, let $M_0$ be the (nonclosed) linear 
span of the monomials in $M$.  $M_0$ is dense in $M$ and 
invariant under the action of all polynomials in 
$S_1,\dots,S_d$.  Thus $M_0$ is a submodule 
of $\mathbb C[z_1,\dots,z_d]\otimes E$ (where 
the latter is considered a finitely generated module 
over the algebra of polynomials $\mathbb C[z_1,\dots,z_d]$).  By 
Hilbert's basis theorem (in the form which asserts that a submodule 
of a finitely generated $\mathbb C[z_1,\dots,z_d]$-module is finitely generated), 
it follows that there is a set of polynomials $f_1,\dots,f_s$ and a 
set $\zeta_1,\dots,\zeta_s\in E$ such that $M_0$ is generated by 
$\{f_1\otimes\zeta_1,\dots,f_s\otimes\zeta_s\}$.  Using the invariance 
of $M_0$ under $\Gamma$ we can 
decompose each $f_j\otimes \zeta_j$ into a finite 
sum of monomials in $M_0$ to obtain the required finite set of generators 
for $M$.  

We conclude that there is a finite set of $d$-tuples 
$\{\nu_1,\dots,\nu_s\}$ in $\mathbb Z^d$,   
satisfying $\nu_j\geq 0$, $\nu_i\neq \nu_j$ for $i\neq j$, 
and a finite 
set of subspaces $E_1,\dots,E_s\subseteq E$ with 
 the property that $M$ is generated as follows
\begin{equation}\label{oldEq2.10}
M=\overline{\text span}
\{f_1(z)z^{\nu_1}\otimes\zeta_1+\dots+f_s(z)z^{\nu_s}\otimes \zeta_s\}, 
%\tag{2.10}
\end{equation}
where $f_1,\dots,f_s$ range over $\mathbb C[z_1,\dots,z_d]$ and 
$\zeta_j\in E_j$.  

We must identify the spectral subspaces of $M$
$$
M(n)=\{\xi\in M: \Gamma(\lambda)\xi=\lambda^n\xi,\quad \lambda\in\mathbb T^d\},
\qquad n\in\mathbb Z^d
$$
in terms of $\nu_1,\dots,\nu_s$ and $E_1,\dots,E_s$.  This requires 
some care since the spaces $E_1,\dots,E_s$ need not be  
mutually orthogonal.  

\begin{lem}\label{oldLem2.11}For every $n\in\mathbb Z^d$ and every 
$k=1,\dots,s$ define $E_k(n)=E_k$ if $\nu_k\leq n$, 
$E_k(n)=\{0\}$ if $\nu_k\nleq n$, and let  
$E(n)=E_1(n)+\dots+E_s(n)$.  Then $M(n)=\{0\}$ except when 
$n\geq 0$, and for $n\geq 0$ we have 
$$
M(n)=z^n\otimes E(n)=\{z^n\otimes\zeta: \zeta\in E(n)\}.  
$$
\end{lem}
\begin{proof}
Clearly $M(n)=\{0\}$ except when $n\geq 0$, so fix $n\geq 0$ 
and let $Q(n)$ be the orthogonal projection 
of $M$ onto $M(n)$.  We 
show that $Q(n)M=z^n\otimes E(n)$.

For each $k=1,\dots,s$ let $M_k\subseteq M$ be the invariant subspace 
generated by $z^{\nu_k}\otimes E_k$,
$$
M_k=\overline{span}\{z^p\otimes E_k: p\geq \nu_k\}.  
$$
If $k$ is such that $ \nu_k\leq n$ then $M_k$ contains 
$z^n\otimes E_k$ and in fact $Q(n)M_k=z^n\otimes E_k$.  
If on the other hand $k$ is 
such that $\nu_k\nleq n$ then 
$n$ cannot belong to the set $\{p\in\mathbb Z^d: p\geq \nu_k\}$ and hence 
$z^n\otimes E$ must be orthogonal to $M_k$.  It follows that 
$Q(n)M_k=\{0\}$ when $\nu_k\nleq n$.  We
conclude that 
$$
Q(n)(M_1+\dots+M_s)=\sum_{k=1}^s z^n\otimes E_k(n)=z^n\otimes E(n).  
$$
After one notes that 
$M_1+\dots+M_s$ must be dense in $M$ because of (\ref{oldEq2.10}), 
the proof is complete.  
\end{proof}

Consider now the subspace $\overline{S_1M}$.  Letting 
$e_1\in\mathbb Z^d$ be the $d$-tuple with components $(1,0,\dots,0)$ we see from 
(\ref{oldEq2.10}) that $\overline{S_1M}$ is generated by 
$$
z^{\nu_1+e_1}\otimes E_1,\dots,  z^{\nu_s+e_1}\otimes E_s, 
$$
and Lemma \ref{oldLem2.11} implies that its spectral subspaces are given 
by 
$$
(\overline{S_1M})(n)=z^n\otimes F(n), \qquad n\geq 0,
$$
where $F(n)$ is defined as the sum $F_1(n)+\dots+F_s(n)$ where 
$F_k(n)=\{0\}$ if $\nu_k+e_k\nleq n$ and 
$F_k(n)=E_k$ if $\nu_k+e_1\leq n$.  

Obviously $F_k(n)\subseteq E_k(n)$ for 
every $k=1,\dots,s$, so that $F(n)\subseteq E(n)$.  
Thus 
$M\ominus S_1M$ decomposes into an orthogonal direct sum 
$$
M\ominus S_1M=\sum_{n\geq 0} E(n)\ominus F(n).  
$$

Finally, let $q_k\in \mathbb Z^+$ be the first component
of $\nu_k=(q_k,*,*,\dots)$, $k=1\dots,s$.  
We claim that any $n\geq 0$ 
for which $E(n)\ominus F(n)\neq \{0\}$ must have its first component 
in the set $\{q_1,\dots,q_s\}$.  Indeed, for such an $n$ 
there must be a $k=1,\dots,s$ such that $F_k(n)\neq E_k(n)$, and 
this implies that $F_k(n)=\{0\}$ and $E_k(n)=E_k$.  The first condition
implies that $\nu_k+e_1\nleq n$ and the second implies that 
$\nu_k\leq n$; hence the first  
component of $n$ must agree with the first component $q_k$ of $\nu_k$.  

Setting $q=\max(q_1,\dots,q_s)$, we conclude that 
$$
M\ominus S_1M\subseteq\overline{\text span}\{z^n\otimes E: n=(n_1,\dots,n_d)\geq 0, 
\quad 0\leq n_1\leq q\},
$$
and Assertion 1 is proved.
\end{proof}

\begin{proof}[proof of Assertion 2]
It is pointed out in (\cite{arvSubalgIII}, Corollary of Proposition 5.3) 
that the operators $S_k$ are 
hyponormal, $S_kS_k^*\leq S_k^*S_k$, $1\leq k\leq d$.  Thus it 
suffices to show that the restriction $C$ of $S_1$ to 
$$
K=\overline{span}\{z^n\otimes\zeta: n_1\leq q,
\quad\zeta\in E\}
$$
satisfies $C^*C\in\mathcal L^p$ for every $p>d$.  $C^*C$ is the 
compression of $S_1^*S_1$ to $K$.  We have seen in the proof 
of Lemma \ref{oldProp2.2} that monomials are eigenvectors for 
$S_1^*S_1$, 
\begin{equation}\label{oldEq2.12}
S_1^*S_1:z^n\otimes \zeta\mapsto 
\frac{n_1+1}{|n|+1}\cdot z^n\otimes \zeta, \qquad n=(n_1,\dots,n_d),
%\tag{2.12}
\end{equation}
where $|n|=n_1+\dots+n_d$.  
Let $\zeta_1,\dots,\zeta_r$ be an orthonormal 
basis for $E$.  Then $K$ is spanned by the orthogonal 
set of all monomials of the form 
$z^n\otimes\zeta_j$ with $n_1\leq q$, $1\leq j\leq r$, and 
from (\ref{oldEq2.12}) it follows that 
\begin{equation}\label{oldEq2.13}
0\leq P_KS_1^*S_1P_K\leq (q+1)P_K((\mathbf 1+N)^{-1}\otimes\mathbf 1_E)P_K,
%\tag{2.13}
\end{equation}
where $N$ is the number operator of $H^2$, the unbounded 
self-adjoint operator 
having the set of monomials as eigenvectors as follows: 
$N: z^n\mapsto |n|z^n$, $n\geq 0$. 
It is known that $(\mathbf 1_{H^2}+N)^{-1}$ belongs to $\mathcal L^p$ 
for every $p>d$ (see formula (5.2) of \cite{arvSubalgIII}).  Since $E$ is 
finite dimensional we have 
$$
(\mathbf 1+N\otimes 1_E)^{-1} = (\mathbf 1_{H^2}+N)^{-1}\otimes \mathbf 1_E
\in\mathcal L^p,\qquad p>d.  
$$
In view of (\ref{oldEq2.13}), 
we conclude that  $P_KS_1^*S_1P_K\in\mathcal L^p$ for every 
$p>d$.
\end{proof}

Lemma \ref{oldProp2.3} follows from Assertions 1 and 2, 
thereby completing the proof of Proposition \ref{oldAss2.1} 
and Theorem \ref{oldThmA}.  
\end{proof}

It remains only to deduce Corollary \ref{oldCor}.  We sketch 
the argument as follows.  Let $(T_1,\dots,T_d)$ be the 
$d$-tuple acting on $M^\perp$ by compression
$$
T_k=(\mathbf 1-P_M)S_k\restriction_{M^\perp}, \quad k=1,\dots,d, 
$$
and let $\mathcal A$ be the $*$-algebra generated by $T_1,\dots,T_d$ and 
the identity operator.  A straightforward argument (that we 
omit) shows that the set of commutators $\{AB-BA: A,B\in\mathcal A\}$ is 
contained in $\mathcal L^p$ iff all the self-commutators 
$T_i^*T_j-T_jT_i^*$, $1\leq i,j\leq d$, belong to $\mathcal L^p$.  

Thus it suffices to show that the self-commutators all belong to 
$\mathcal L^p$ for $p>d$.  Writing $P$ for the projection onto $M$, 
$P^\perp$ for $\mathbf 1-P$ and $T_k=P^\perp S_k P^\perp$ we have 
\begin{align}
T_i^*T_j&=P^\perp S_i^*(\mathbf 1-P)S_jP^\perp =
P^\perp S_i^*S_jP^\perp - P^\perp S_i^*PS_jP^\perp\\
&=P^\perp S_i^*S_jP^\perp - A_iA_j^*, 
\end{align}
where $A_i=P^\perp S_i^*P$, and 
$
T_jT_i^* =P^\perp S_jP^\perp S_i^*P^\perp = P^\perp S_jS_i^*P^\perp.  
$
Thus 
$$
[T_i^*,T_j]=P^\perp [S_i^*,S_j]P^\perp - A_iA_j^*.  
$$
Proposition \ref{oldAss2.1} implies that $A_i\in\mathcal L^p$ for $p>2d$, and 
hence $A_iA_j^*\in\mathcal L^p$ for $p>d$.  Finally, 
according to Proposition 5.3 of \cite{arvSubalgIII} we have
$[S_i^*,S_j]\in\mathcal L^p$
for $p>d$, and the desired conclusion follows.  

\section{Submodules and Quotients}\label{S:subM}

Every $d$-contraction $\bar A=(A_1,\dots,A_d)$ has 
a defect operator 
$$
\Delta_{\bar A}=\mathbf 1-(A_1A_1^*+\dots+A_dA_d^*), 
$$
and one has $0\leq \Delta_{\bar A}\leq \mathbf 1$.  
While this notation differs from that of \cite{arvSubalgIII} 
where $\Delta_{\bar A}$ was defined as the square root of 
$\mathbf 1-(A_1A_1^*+\dots+A_dA_d^*)$, 
it is better suited for our purposes here.  
We use the 
traditional notation $[X,Y]$ to denote the 
commutator $XY-YX$ of two operators $X,Y$. 

Given an invariant subspace $M\subseteq H$ for a 
$d$-contraction $\bar A$, the restriction of $\bar A$ 
to $M$ and the compression of $\bar A$ to the quotient 
$H/M$ define two new $d$-contractions.  In this section 
we examine the relationships between these three 
multioperators.  We identify the quotient Hilbert 
space $H/M$ with the orthocomplement $M^\perp$ 
of $M$ in $H$, and its associated $d$-contraction with 
the $d$-tuple obtained by compressing 
$(A_1,\dots,A_d)$ to $M^\perp$.   

\begin{prop}\label{prop4}
Let $\bar A=(A_1,\dots,A_d)$ be a $d$-contraction acting on 
a Hilbert space $H$, 
let $M$ be a closed 
$\bar A$-invariant subspace with projection $P:H\to M$,  
and let $\bar B=(B_1,\dots,B_d)$ and 
$\bar C=(C_1,\dots,C_d)$ be the $d$-contractions obtained, 
respectively, by 
restricting $\bar A$ to $M$ and compressing $\bar A$ to 
$M^\perp$.  
Writing $P^\perp$ for the projection onto 
the subspace $M^\perp\subseteq H$, we have the 
following formulas relating various commutators and 
the three defect 
operators $\Delta_{\bar A}$, $\Delta_{\bar B}$, $\Delta_{\bar C}$.  
\begin{align}
[B_j,B_k^*] P=&-[P,A_j][P,A_k]^*+P[A_j,A_k^*]P\label{Eq4}\\
[C_j,C_k^*] P^\perp=
&[P,A_k]^*[P,A_j] + P^\perp[A_j,A_k^*]P^\perp\label{Eq5}\\
\Delta_{\bar B} P=&
P\Delta_{\bar A}P+ \sum_{k=1}^d[P,A_k][P,A_k]^*\label{Eq6}\\
\Delta_{\bar C} P^\perp=&
P^\perp\Delta_{\bar A}P^\perp\label{Eq7}.  
\end{align}
\end{prop}

\begin{proof}
To verify (\ref{Eq4}), we write 
\begin{align*}
[B_j,B_k^*]P=&A_jPA_k^*P-PA_k^*A_jP = 
A_jPA_k^*P-PA_jA_k^*P +P[A_j,A_k^*]P\\
=&-PA_jP^\perp A_k^*P+P[A_j,A_k^*]P.  
\end{align*}
Since $PA_jP^\perp =PA_j-A_jP$, we have 
$PA_jP^\perp A_k^*P=[P,A_j][P,A_k]^*$, 
and (\ref{Eq4}) follows.  

(\ref{Eq5}) follows similarly, 
after using $P^\perp A_j P^\perp = P^\perp A_j$ 
to write 
\begin{align*}
[C_j,C_k^*]P^\perp 
&=P^\perp A_jA_k^*P^\perp-P^\perp A_k^*P^\perp A_jP^\perp 
\\
&=P^\perp A_k^*A_jP^\perp -P^\perp 
A_k^*P^\perp A_jP^\perp +P^\perp [A_j,A_k^*]P^\perp\\
&=P^\perp A_k^*PA_jP^\perp + P^\perp [A_j,A_k^*]P^\perp.  
\end{align*}
(\ref{Eq5}) follows after one notes that 
$P^\perp A_k^*PA_jP^\perp=[P,A_k]^*[P,A_j]$.  

To prove (\ref{Eq6}), one writes $\Delta_{\bar B}P$ as follows,
\begin{align*}
P-\sum_{k=1}^dA_kPA_k^*&
=P\Delta_{\bar A}P+\sum_{k=1}^d PA_k(\mathbf 1-P) A_k^*P\\
&
=P\Delta_{\bar A}P+\sum_{k=1}^d[P,A_k][P,A_k]^*, 
\end{align*}
and (\ref{Eq7}) follows similarly, since 
$$
P^\perp -\sum_{k=1}^dP^\perp A_k P^\perp A_k^*P^\perp =
P^\perp \Delta_{\bar A}P^\perp +\sum_{k=1}^dP^\perp A_kPA_k^*P^\perp
=P^\perp\Delta_{\bar A}P^\perp.  
$$
That completes the proof.
\end{proof}

\begin{cor}Let $\bar A$, $\bar B$, $\bar C$ satisfy 
the hypotheses of Proposition \ref{prop4}.  
Then for every $p$ satisfying $1\leq p\leq \infty$, the 
following are equivalent:
\begin{enumerate}
\item[(i)]Both defect operators $\Delta_{\bar B}$ and 
$\Delta_{\bar C}$ belong to $\mathcal L^p$.  
\item[(ii)] $\Delta_{\bar A}$ belongs 
to $\mathcal L^p$ and $[P_M,A_k]\in\mathcal L^{2p}$.  
$1\leq k\leq d$.  
\end{enumerate}
\end{cor}

\begin{proof}The implication (ii)$\implies$(i) is an 
immediate consequence of the formulas (\ref{Eq6}) and (\ref{Eq7}).  

(i)$\implies$(ii):  We write $P$ for $P_M$.  
From (\ref{Eq6}) and (\ref{Eq7}), together 
with the fact that the right side of (\ref{Eq6}) is a sum of 
positive operators, we may conclude that all of the 
operators 
$$
P\Delta_{\bar A}P,\   P^\perp\Delta_{\bar A}P^\perp, \ 
[P,A_1][P,A_1]^*,\dots,[P,A_d][P,A_d]^*
$$ 
belong to $\mathcal L^p$.  
By (\ref{Eq12}) we have 
$[P,A_k]\in\mathcal L^{2p}$, $1\leq k\leq d$.  
Another application of (\ref{Eq12}) 
shows that both $\sqrt{\Delta_{\bar A}}P$ 
and $\sqrt{\Delta_{\bar A}}P^\perp$ belong to $\mathcal L^{2p}$.  
The latter two operators sum to 
$\sqrt{\Delta_{\bar A}}\in\mathcal L^{2p}$, 
and therefore $\Delta_{\bar A}\in\mathcal L^p$.  
\end{proof}

We now apply Proposition \ref{prop4} to obtain concrete information 
about the examples of greatest interest for us, namely the 
submodules and quotients that are associated with 
pure finite rank $d$-contractions.  

\begin{thm}\label{thm2}
Let $M\subseteq H^2\otimes \mathbb C^r$ 
be an invariant subspace of the $d$-shift 
$\bar S=(S_1,\dots,S_d)$ 
of finite rank $r$ and let $\bar B$ and $\bar C$ be, 
respectively, the restriction of $\bar S$ to $M$ 
and the compression of $\bar S$ to $M^\perp$.   
Then for every $p$ satisfying 
$d<p\leq \infty$, the following are equivalent:   
\begin{enumerate}
\item[(i)]The defect operator of $\bar B$ belongs to $\mathcal L^p$.  
\item[(ii)]$[B_j,B_k^*]\in \mathcal L^p$, $1\leq j,k\leq d$.  
\item[(iii)]$[C_j,C_k^*]\in\mathcal L^p$, $1\leq j,k\leq d$.  
\item[(iv)]$[P_M,S_k]\in\mathcal L^{2p}$, $1\leq k\leq d$.  
\end{enumerate}
If (i)--(iv) are satisfied for some $p\in(d,\infty]$, then 
both $\bar B$ and $\bar C$ are Fredholm multioperators, and 
the indices of their Dirac operators are related by 
\begin{equation}\label{Eq11}
\ind D_{\bar B\,+}+\ind D_{\bar C\,+}=(-1)^d \cdot r.  
\end{equation}
\end{thm}

\begin{proof}
It was shown in \cite{arvSubalgIII} that the self-commutators 
$[S_j,S_k^*]$ belong to $\mathcal L^q$ for every $q>d$; and 
of course, the defect operator of $\bar S$ is a projection 
of rank $r$, belonging to $\mathcal L^q$ for 
every $q\geq1$.  With these observations in hand, one sees  
from (\ref{Eq6}) that $\Delta_{\bar B}\in\mathcal L^p$ if 
and only if 
\begin{equation}
\sum_{k=1}^d[P,S_k][P,S_k]^*\in\mathcal L^p\label{Eq10}, 
\end{equation}
where $P$ denotes $P_M$.   
Similarly, (\ref{Eq4}) and (\ref{Eq5}) 
show that the assertions 
(ii) and (iii) are equivalent, respectively, 
to the 
assertions 
\begin{align}
[P,S_j][P,S_k]^*&\in\mathcal L^p,\qquad 1\leq j,k\leq d\label{Eq8}\\
[P,S_k]^*[P,S_j]&\in\mathcal L^p, \qquad 1\leq j,k\leq d\label{Eq9}. 
\end{align}
Thus, the problem of showing that (i)--(iv) are equivalent is 
reduced to that of showing that each of the assertions 
(\ref{Eq10}), (\ref{Eq8}) and (\ref{Eq9}) is equivalent to the 
assertion $[P,S_k]\in\mathcal L^{2p}$, $1\leq k\leq d$.  
That 
is a straightforward consequence of the 
elementary equivalences 
(\ref{Eq12}).  

To sketch the proof of (\ref{Eq11}), note first that Proposition 
\ref{prop1}, together with (ii), (iii), and the known 
essential normality of the $d$-shift, imply 
that the three Dirac operators 
$D_{\bar S}$, $D_{\bar B}$ and $D_{\bar C}$ are 
Fredholm.  By property (iv), the commutators $PS_k-S_kP$, 
$1\leq k\leq d$, are compact.  It follows that the $d$-tuple 
$\bar S$ is unitarily equivalent to a compact perturbation of 
the direct sum of $d$-tuples 
$\bar B\oplus \bar C$.   
In turn, this implies that 
$D_{\bar S\,+}$ is unitarily equivalent to 
a compact perturbation of the direct sum 
of Fredholm operators 
$D_{\bar B\,+}\oplus D_{\bar C\,+}$.  By stability of 
the Fredholm index under compact perturbations, we have 
$$
\ind D_{\bar B\,+}+\ind D_{\bar C\,+}=\ind D_{\bar S\,+}.     
$$
Finally, from Theorem B of \cite{arvDirac} that relates 
the index of a finite rank graded $d$-contraction to 
its curvature invariant, we can compute the right 
side of the preceding formula
$$
\ind D_{\bar S\,+}=(-1)^dK(\bar S).    
$$
Since the curvature of a finite direct sum of copies 
of the $d$-shift is known to be its rank \cite{arvCurv}, 
formula (\ref{Eq11}) follows.  \end{proof}

\section{Concluding Remarks and Conjectures}\label{S:conc}

We expect that some variation of Corollary 
\ref{oldCor} should  
hold under much more general circumstances, 
and we now discuss these issues.

\vskip0.1in
{\bf Conjecture A.}
{\em 
Let $M$ be 
a closed invariant subspace for the $d$-shift 
$\bar S=(S_1,\dots,S_d)$ of finite 
rank $r$, acting on 
$H^2\otimes \mathbb C^r$.  Assume that $M$ is generated 
by a  set of 
vector polynomials 
in $\mathbb C[z_1,\dots,z_d]\otimes \mathbb C^r$,  
each of which is homogeneous of some degree.  
Then $P_MS_k-S_kP_M$ 
is compact for every
$k=1,\dots,d$.  
}
\vskip0.1in
Note that because of Hilbert's basis theorem, one 
may assume that $M$ is generated by a {\it finite} set 
of homogeneous vector polynomials .

By Theorem \ref{thm2}, Conjecture A implies that the pure 
$d$-contraction $\bar T=(T_1,\dots,T_d)$ obtained by compressing 
$\bar S$ to $M^\perp$ is a Fredholm multioperator, 
and as we have seen in Section \ref{S:intro}, the index formula 
(\ref{Eq2}) implies that the curvature invariant $K(\bar T)$ 
is  stable in these cases.  The space $X$ 
appearing in the exact sequence of \cstar s 
$$
0\longrightarrow \mathcal K\longrightarrow 
C^*(T_1,\dots,T_d)+\mathcal K
\longrightarrow 
C(X)\longrightarrow 0
$$
would now be associated with a nontrivial 
algebraic set in projective space.    

\begin{rem}[Evidence for Conjecture A]\label{rem17}
Theorem \ref{oldThmA} 
implies that Conjecture A
is true when the homogeneous polynomials are 
monomials.  
Moreover, Conjecture A is true in 
two dimensions.  Indeed, a recent result of Kunyu Guo 
(Theorem 2.4 of \cite{guoDef}) implies that, 
in the context of Conjecture A for dimension $d=2$, 
the $2$-contraction obtained by restricting 
$(S_1,S_2)$ to $M$ has the property that its 
defect operator belongs to $\mathcal L^p$ for every 
$p>1$.   By Theorem \ref{thm2}, 
Conjecture A is true when $d=2$.  

Finally, there are a few other 
classes of (unpublished) examples in 
arbitrary dimension $d$ involving homogeneous 
polynomials for which one can decide the 
issue, and these too support Conjecture A.  
\end{rem}

Stephen 
Parrott has shown  \cite{parrottCurv} 
that a pure finite rank single contraction  
is a Fredholm operator and 
(\ref{Eq1}) holds; R. N. Levy improved this 
in \cite{levyRN}.  However, some of 
the examples that occur in this 
one-dimensional setting are not essentially 
normal.  Thus, one cannot expect the 
conclusion of Conjecture A to hold 
for arbitrary invariant subspaces 
$M\subseteq H^2\otimes \mathbb C^r$ in higher dimensional 
cases $d>1$.  However, we believe that 
the following two ``ungraded" relatives of 
Conjecture A are well-founded.  
\vskip0.05in
{\bf Conjecture B.}
{\em 
Let $M\subseteq H^2\otimes \mathbb C^r$ 
be an invariant subspace of the 
$d$-shift of finite rank $r$ that is generated 
by a set of vector polynomials.  Then 
$P_MS_k-S_kP_M$ is compact for $1\leq k\leq d$.  
}
\vskip0.05in
Assuming the result of Conjecture B, one obtains 
a Fredholm multioperator by compressing the 
$d$-shift to $M^\perp$.  In such cases, one  
would expect the following conjecture to hold; 
the result would generalize 
the index formula (\ref{Eq1}) to the ungraded case.   
\vskip0.05in
{\bf Conjecture C.}
{\em 
The index formula (\ref{Eq1}) holds for the 
finite rank $d$-contraction $\bar T$ obtained by compressing 
the $d$-shift of rank $r$ to $M^\perp$, whenever $M$ is generated 
by vector polynomials and $\bar T$ is essentially normal.  
}
\vskip0.05in
More generally, it is natural to ask if every finite rank pure 
$d$-contraction is Fredholm and satisfies 
the index formula (\ref{Eq1}).  While there 
is scant evidence 
to illuminate  
these questions in general, 
Parrott's work \cite{parrottCurv} implies that both answers 
are yes in the one-dimensional cases $d=1$:  

\vskip0.05in
{\bf Problem D.}
{\em 
Let $\bar T$ be a finite rank $d$-contraction.  Is 
$\bar T$ is a Fredholm $d$-tuple?  Does the index 
formula (\ref{Eq1}) hold?  
}
\vskip0.05in
We expect that significant progress on Problem D will 
require further development of the 
theory of Fredholm multioperators.

%\vfill

%\bibliography{bibData}  %Remove this line when finished, see below.  
\bibliographystyle{alpha}

\newcommand{\noopsort}[1]{} \newcommand{\printfirst}[2]{#1}
  \newcommand{\singleletter}[1]{#1} \newcommand{\switchargs}[2]{#2#1}

%\bibliography{schurHorn}				%%Insert this when finished with ms. 
\end{document}